\documentclass[1p,12pt]{elsarticle}

\usepackage{amssymb}
\usepackage{amsmath}
\usepackage{graphicx}

\journal{XXXXXXX}
\begin{document}

\begin{frontmatter}

\title{Locality preserving homogeneous Hilbert curves by use of arbitrary kernels}

\author[eer]{C. P\'erez-Demydenko}

\author[baf]{I. Brito Reyes}

\author[phys]{E. Estevez-Rams\corref{cor1}}
\ead{estevez@imre.oc.uh.cu}

\author[baf]{B. Arag\'on Fern\'andez}

\address[eer]{Instituto de Ciencias y Tecnolog\'ia de Materiales (IMRE), University of Havana , San Lazaro y L. CP 10400. La Habana. Cuba.}

\address[baf]{Universidad de las Ciencias Inform\'aticas (UCI), Carretera a San Antonio, Boyeros. La Habana. Cuba.}

\address[phys]{Physics Faculty-Instituto de Ciencias y Tecnolog\'ia de Materiales (IMRE), University of Havana, San Lazaro y L. CP 10400. La Habana. Cuba.}

\cortext[cor1]{Corresponding author at: Instituto de Ciencias y Tecnolog\'ia de Materiales, University of Havana (IMRE), San Lazaro y L. CP 10400. La Habana. Cuba.}

\begin{abstract}
Homogeneous Hilbert curves (HHC) in two dimensions are generalized by introducing the construction of the space filling curves from the same affine transformations but using an arbitrary kernel, we call such curves HHCK. The new curves are still space filling that comply with the nesting condition but violates the adjacency property. The freedom of building new HHC curves with arbitrary kernels, is only limited by the constrain that the chosen kernel must allow a well behaved connectivity between quadrants. Two examples of such curves are discussed. 

The important property of locality preservation in space filling curve mapping is discussed. Besides the common used dilation factor, the paper introduces and discuss difference map as a site locality measure, allowing to describe locality preservation in a more detailed way than dilation factors. The strength of such analysis is proven and global descriptors from the difference map are derived. Locality of all HHC curves and the two discussed HHCK curves are studied by difference map.  
\end{abstract}



\date{\today}
\end{frontmatter}

\section{Introduction}

Space filling curves allow a surjective mapping between the one-dimensional space and the $d$-dimensional one: $\mathbb{R}\longrightarrow\mathbb{R}^{d}$. In the two-dimensional case, one of the most studied space filling curve has been the Hilbert curve. Hilbert curve allows to construct an onto, but not one-to-one, mapping between the unit segment $I=\{ t | 0 \leqslant t \leqslant 1\}$ and the unit square $Q=\{ (x,y) | 0 \leqslant x \leqslant 1,0 \leqslant y \leqslant 1\}$ \cite{sagan94}. This mapping, as any other space filling curve mapping, is not a bijective one, as proved by Netto \cite{sagan94}, a fact that can be seen as result of the infinite application of the same affine transformation over an initial curve. Nevertheless, their $n$th-order iterations or $n$th approximations, may indeed be seen as a bijective mapping between sub-segments of $I$ and  sub-squares of $Q$ to an arbitrary fine granularity with increasing $n$ \cite{bauman06} . This, together with the good locality  preservation through all the iterations, has opened a wealth of applications of the Hilbert curve in computer sciences \cite{chen05,chen11,bially69}, image processing \cite{songa02,liang08}, data representation \cite{anders09} and other fields.

The Hilbert curve is uniquely defined using the so called nesting and adjacency conditions, and restricting the first and last subintervals of $I$ to correspond to the lower left and right, respectively, sub-squares of $Q$ \cite{sagan94}. Changing the initial and ending point, Moore \cite{moore00}  introduced a new curve that now bears his name. More recently, Liu \cite{liu04} described four new curves,  that introduces a new approach for constructing other mappings between sub-segments of $I$ and  sub-squares of $Q$. The affine transformations involved to obtain a $n$-order curve, involves only one iteration over the $(n-1)$-order Hilbert curve. 

P\'erez-Davidenko et al. \cite{estevez_may13}, have further developed the idea by introducing the concept of homogeneous Hilbert Curves (HHC) in two dimensions which can be proper or improper. Homogeneity implies that only one set of rules are applied to the $n$-order Hilbert curve, in order to generate the   $(n+1)$-order curve. The development is based on the idea that $n$-order space filling curves can be built from $(n-1)$-order Liu curves, if the latter shows the correct quadrant connectivity. It was found that only the fourth Liu curve can be used for such purpose, and six additional curves can be added to complete twelve HHC in two dimensions \cite{estevez_may13}. 

The importance of finding new space filling curves is due to the fact that several applications can benefit of a broader set  from where it can draw candidate solutions. An example of such application is the traveling salesman problem by use of space filling curves \cite{gao94,bartholdi82}.  

In this contribution we report two results. First, we show how Hilbert curve based space filling curves can be further generalized by relaxing the condition that they are constructed from the iteration  of the first order Hilbert curve. Indeed, by defining HHC from the successive iteration of an arbitrary \textit{kernel}, new space filling curves can be built, which still serve as mappings between sub-segments of $I$ and  sub-squares of $Q$, but which may additionally violate the adjacency condition. The locality preservation properties of such curves will depend on the chosen kernel and the affine transformations involved. The freedom of building new HHC curves with arbitrary kernels, is only hindered by the constrain that the chosen kernel must allow a well behaved connectivity between quadrants, in the same manner as already explained in  \cite{estevez_may13}. Two such kernels will be discussed as examples. 

The second result is to explore the locality properties of all the presented curves. We evaluate locality by the dilation factor, closely related to the H\"older condition. Dilation factor has been used to analyze locality in the infinite limit of space filling curves \cite{shchepin08}. As the dilation factor is a global descriptor of locality, we introduce  the difference map, as a site measure of departure from locality in space filling curves. 

\section{Homogeneous Hilbert Curves}\label{sec:HHC}

For the purpose of this paper, a 2D space filling curve is a surjective mapping of the unit interval $I$ onto the unit square $Q$, $I \longrightarrow Q$. In the Hilbert construction, we start by partitioning the unit interval into four non-overlapping subintervals, and the unit square into four congruent non-overlapping subsquares. A mapping is made between each subinterval $I_{i}$ of $I$ and each subsquare $Q_{j}$ of $Q$. For each subinterval and its corresponding subsquare, the procedure can be repeated recursively, treating each of them as a new interval and a new square. While repeating this procedure, continuity must be preserved by assuring that adjacent subsquares corresponds to adjacent subintervals. A Hilbert curve of order $n$ will correspond to the partition of the unit interval and the unit square into $4^n$ subintervals and subsquares, respectively. The Hilbert curve of order $n$ is uniquely determined by giving the initial and final subsquares corresponding to the initial and final subintervals.  

In \cite{estevez_may13} it was shown, that there exist 12 HHC, half of them were called proper and the other half, improper. The HCC have been labeled with an index $\nu$ going from 0 to 11. The values of $\nu=$0, 1, 2, 3, 4 and 5 correspond to the proper HHC, and the rest of the values, to the improper ones. The value $\nu=0$ correspond to the classic Hilbert curve, $\nu$=1 to the Moore curve, and the values $\nu$=2,3,4 and 5 to the four Liu curves. For these curves, we may build infinitely approximations, getting for increasing iteration, a finer granularity in the spanning of the unit square.  We may use the order $n$ as a label, thus $_{\nu}H_{n}$ will represent the $\nu$ HHC of order $n$. In general, the order of the curve can be larger or equal to the number of iterations performed in the building algorithm.

Proper HHC are generated iteratively using specific affine transformations $_{\nu}q_{i}$ that maps one subinterval to the $i$th quadrant. The affine transformations operators have a rotation part $\mathbf{U}$, given by a $2\times2$ integer matrix, and a translation factor, given by a vector $\mathbf{t}$, scaling occurs by $1/2$. The general form of the operator is thus of the form 
\begin{equation*}
q\left (\begin{array}{l}x \\ y\end{array}\right )=\frac{1}{2}\mathbf{U}\cdot \left (\begin{array}{l}x \\ y\end{array}\right )+\frac{1}{2}\mathbf{t},
\end{equation*}
more details can be found in \cite{estevez_may13}, where the complete set of affine transformations is reported. For completeness, Table \ref{tbl:affine} gives the $\mathbf{U}$ and $\mathbf{t}$ for all possible affine transformations.

Quadrants in each subsquare are numbered in clockwise sense, starting by the right lower corner. For the $n$-order approximation of the $_{\nu}H$ curve, the first $n-1$ iterations will use the affine transformation of the Hilbert curve $_{0}H$, but the last affine transformation to apply, will be in general, a different one. That is, for the proper HHC, the following will hold:

\begin{equation}
_{\nu}H_{n} =_{\nu}H\otimes\;_{0}H_{n-1},
\end{equation}
where the symbol $\otimes$ means \textquotedblleft acting over the curve ...\textquotedblright, or what is the same, transforming the points $(x,y)$ which define it. It can be seen that for $\nu$=0 we simply get the classic Hilbert curve of order $n$. By construction, $_{\nu}H_{1}\equiv_{0}H_{1}$.

\begin{figure}
\includegraphics[scale=0.75]{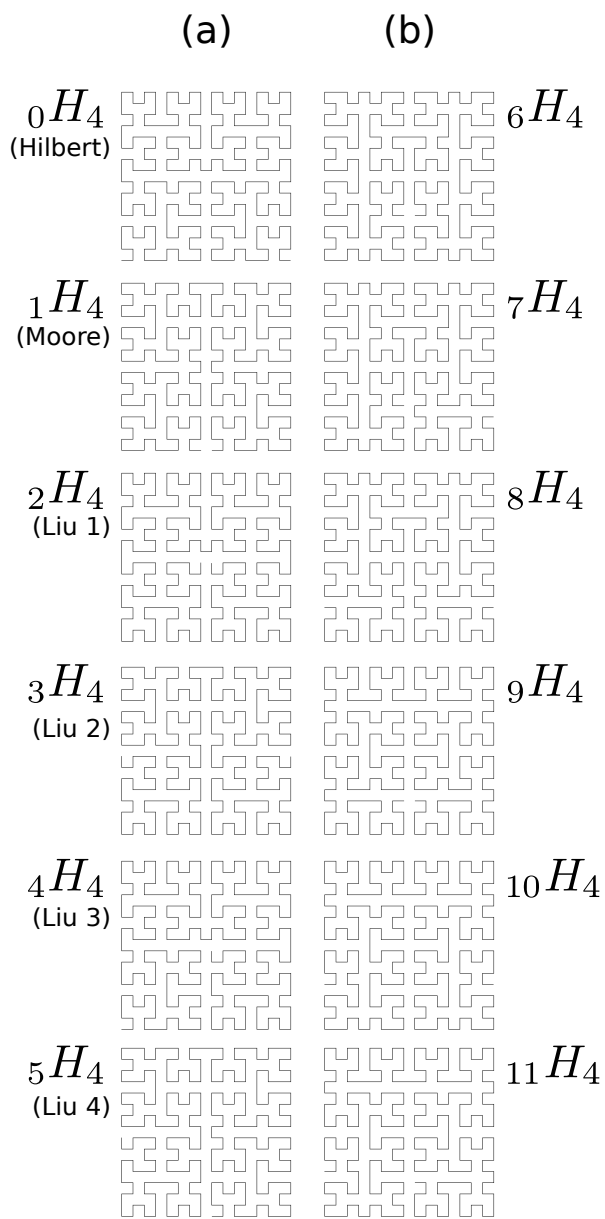}
\caption{Homogeneous Hilbert curves of order 4. (a) proper, (b) improper.
}\label{fig:HHC4}
\end{figure}

Figure \ref{fig:HHC4}a shows the proper Hilbert curves of order $n=4$.

On the other hand, the $_{\nu}H_{n}$ ($n=6,7,\dots11$) improper HHC, are generated by applying the affine transformations $_{\nu}q$ to the proper curve $_{5}H_{n-1}$, i.e. over the fourth Liu curve of order $(n-1)$:

\begin{equation}
_{\nu}H_{n}  =_{\nu}H\otimes\;_{5}H_{n-1} 
\end{equation}

Figure \ref{fig:HHC4}b shows the improper Hilbert curves of order $n=4$.

Such construction is only possible if, besides the affine mappings, a new operation is introduced, the reversion operation, that exchanges the entry and exit point of a curve (detailed description of the full construction algorithm can be found in \cite{estevez_may13}). Since $_{5}H_{n-1} =\;_5H\otimes\;_{0}H_{n-2}$, we will now have for $\nu>5$
\begin{equation*}
_{\nu}H_{1}\equiv_{0}H_{1} 
\end{equation*}
\begin{equation*}
_{\nu}H_{2}\equiv_{5}H_{2}. 
\end{equation*}

\section{Hilbert Curves with arbitrary kernel}\label{sec:kernels}

\begin{figure}
\includegraphics[scale=1.2]{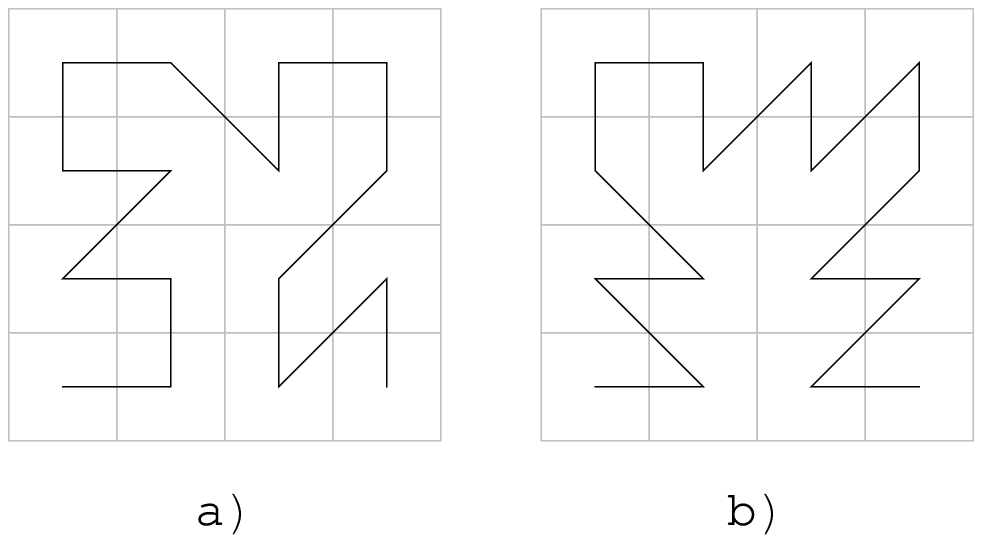}
\caption{a)The mouse and b) frog kernel.
}\label{fig:exampleKernels}
\end{figure}

As described in the previous section, all HHC are generated by successive iterations over the Hilbert curve of order $1$. We  may thus think of $_{0}H_{1}$ as a kernel upon which the other curves are constructed. It can be easily realized that any well connected curve that start at the lower left quadrant and end at the lower right one, will also serve as a kernel. For the first partition of the unit square, there exists only one of such curves, $_{0}H_{1}$, but for higher order partitions there may be several choices. We may, for instance, take as kernels the second order curves showed in Figure \ref{fig:exampleKernels}, and construct with them, twelve different space filling curves by applying the same construction that generated the HHC. We will denote the curves constructed by using an arbitrary kernel by HHCK.  

In general, when using arbitrary kernels, the adjacency condition may be disobeyed, since adjacent sub-intervals in the unit interval $I$ might correspond not to edge, but corner sharing sub-squares of the unit square $Q$. Still, the obtained curve of order $n$ will define a mapping between the sub-intervals and the sub-squares, and therefore may be readily applied for the same tasks as the HHC.

\begin{figure}
\includegraphics[scale=0.8]{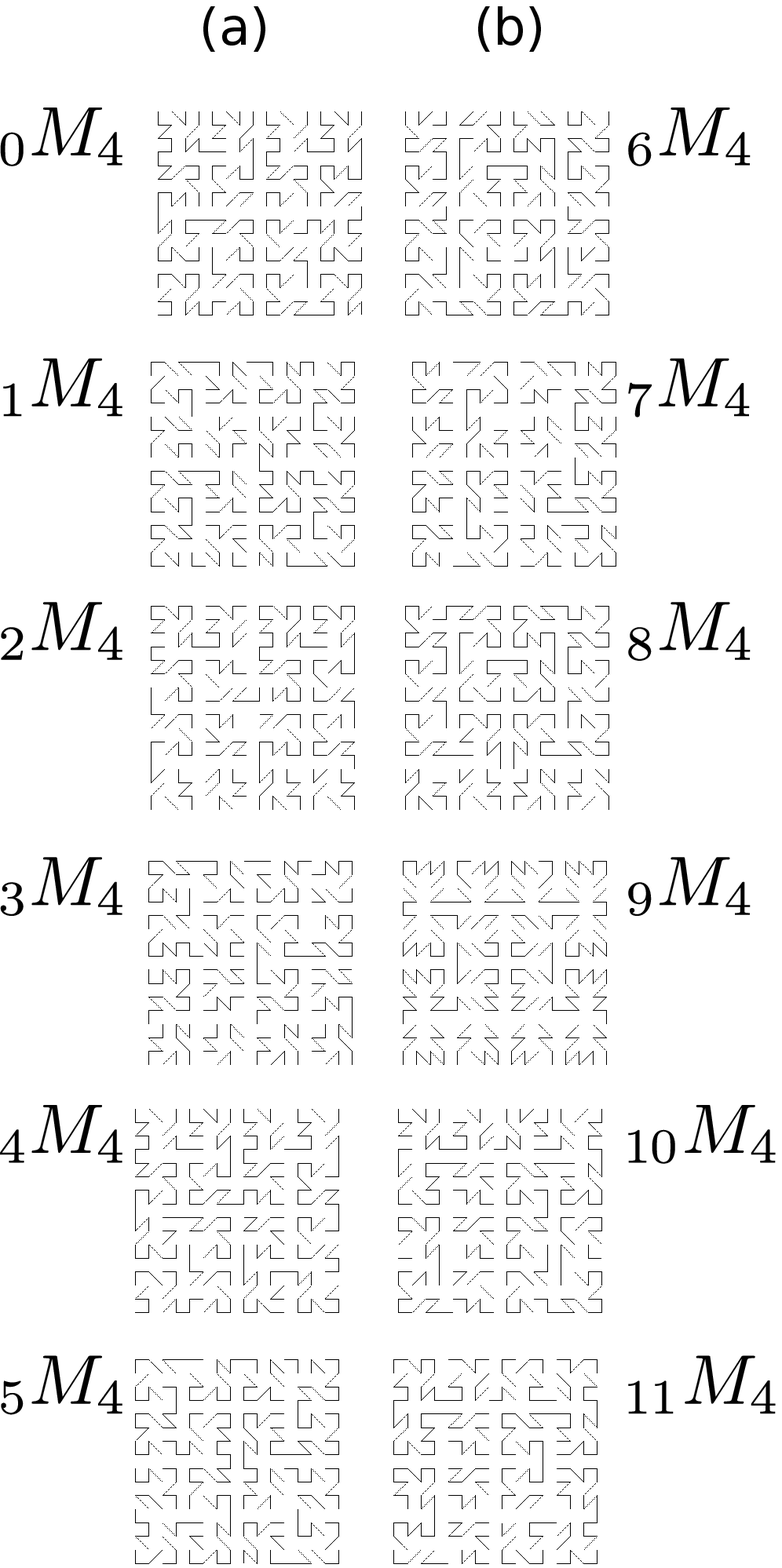}
\caption{The 12 HHCK of order 4 built with the mouse kernel ( see figure \ref{fig:exampleKernels}).a) Proper curves, b) improper curves.
}\label{fig:HHCKmouse}
\end{figure}

\begin{figure}
\includegraphics[scale=0.8]{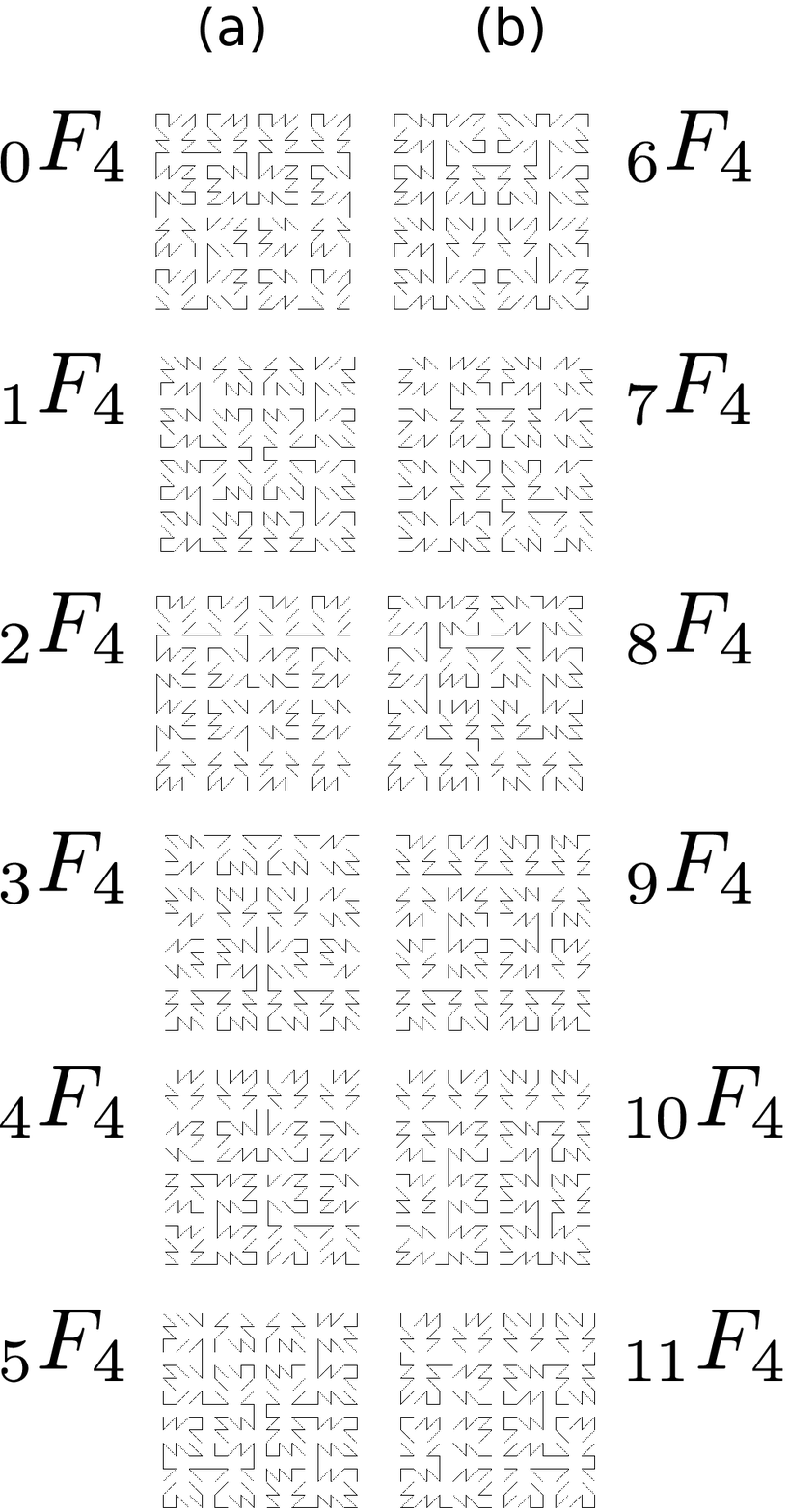}
\caption{The 12 HHCK of order 4 built with the frog kernel ( see figure \ref{fig:exampleKernels}).a) Proper curves, b) improper curves
}\label{fig:HHCKfrog}
\end{figure}

We called the two kernels shown in Figure \ref{fig:exampleKernels}a and b,  the \textit{mouse kernel} and the \textit{frog kernel} and the symbol used for the crresponding curves will be ''M'' and ''F'', respectively. Figures \ref{fig:HHCKmouse} and \ref{fig:HHCKfrog} shows the order four curves built from the mouse and frog kernel respectively.

\begin{figure}
\includegraphics[scale=.3]{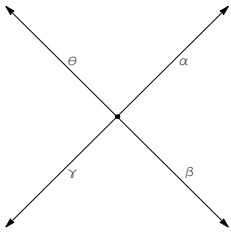}
\caption{The four new diagonal strokes for arbitrary kernel HHC.
}\label{fig:diagonalStrokes}
\end{figure}

In \cite{estevez_may13} a tag system was described for the 12 HHC, a similar description can be devised for the new curves based of the mouse and frog kernel. In order to do so, we must take into account that curves are being constructed with an arbitrary kernel, and therefore extend the definitions of the used morphisms to include the  diagonal strokes that arbitrary kernels introduce. There exist four diagonal strokes that we label $\alpha$, $\beta$, $\gamma$ an $\theta$, for the right-up, right-down, left-down and right-up strokes, respectively (Figure \ref{fig:diagonalStrokes}). The required complete set of stroke transforming operators will be
\begin{equation}
 \begin{array}{llll}
  \delta_{o}(u)=r & \delta_{o}(r)=u & \delta_{o}(d)=l & \delta_{o}(l)=d \\
  \delta_{o}(\alpha)=\alpha & \delta_{o}(\beta)=\theta & \delta_{o}(\gamma)=\gamma & \delta_{o}(\theta)=\beta \\
  \\
  \delta_{a}(u)=l & \delta_{a}(r)=d & \delta_{a}(d)=r & \delta_{a}(l)=u \\
  \delta_{a}(\alpha)=\gamma & \delta_{a}(\beta)=\beta & \delta_{a}(\gamma)=\alpha & \delta_{a}(\theta)=\theta \\
  \\
  \delta_{g}(u)=l & \delta_{g}(r)=u & \delta_{g}(d)=r & \delta_{g}(l)=d \\
  \delta_{g}(\alpha)=\theta & \delta_{g}(\beta)=\alpha & \delta_{g}(\gamma)=\beta & \delta_{g}(\theta)=\gamma \\
  \\
  \delta_{x}(u)=r & \delta_{x}(r)=d & \delta_{x}(d)=l & \delta_{x}(l)=u \\
  \delta_{x}(\alpha)=\beta & \delta_{x}(\beta)=\gamma & \delta_{x}(\gamma)=\theta & \delta_{x}(\theta)=\alpha \\
  \\
  \delta_{f}(u)=d & \delta_{f}(r)=l & \delta_{f}(d)=u & \delta_{f}(l)=r \\
  \delta_{f}(\alpha)=\gamma & \delta_{f}(\beta)=\theta & \delta_{f}(\gamma)=\alpha & \delta_{f}(\theta)=\beta \\
  \\
  \delta_{m}(u)=d & \delta_{m}(r)=r & \delta_{m}(d)=u & \delta_{m}(l)=l \\
  \delta_{m}(\alpha)=\beta & \delta_{m}(\beta)=\alpha & \delta_{m}(\gamma)=\theta & \delta_{m}(\theta)=\gamma \\
  \\
  \delta_{y}(u)=u & \delta_{y}(r)=l & \delta_{y}(d)=d & \delta_{y}(l)=r \\
  \delta_{y}(\alpha)=\theta & \delta_{y}(\beta)=\gamma & \delta_{y}(\gamma)=\beta & \delta_{y}(\theta)=\alpha \\
 \end{array}
\end{equation}

and the tag system for each curve will be the same as the ones reported in \cite{estevez_may13}:

\begin{equation}
\begin{array}{l}
_{0}h_{n+1}=\delta_{o} (_{0}h_{n}) \; u \; _{0}h_{n} \; r \; _{0}h_{n} \; d \; \delta_{a} (_{0}h_{n})    \\
_{1}h_{n+1}=\delta_{g}(_{0}h_{n})  \; u \; \delta_{g}(_{0}h_{n})  \; r \; \delta_{x}(_{0}h_{n}) \; d \; \delta_{x}(_{0}h_{n}) \\
_{2}h_{n+1}=\delta_{f} (_{0}h_{n}) \; u \; _{0}h_{n}              \; r \; _{0}h_{n}             \; d \; \delta_{f} (_{0}h_{n}) \\
_{3}h_{n+1}=\delta_{m}(_{0}h_{n})  \; u \; \delta_{g}(_{0}h_{n})  \; r \; \delta_{x}(_{0}h_{n}) \; d \; \delta_{m}(_{0}h_{n}) \\
_{4}h_{n+1}=\delta_{o} (_{0}h_{n}) \; u \; _{0}h_{n}              \; r \; _{0}h_{n}             \; d \; \delta_{f} (_{0}h_{n}) \\
_{5}h_{n+1}=\delta_{m}(_{0}h_{n})  \; u \; \delta_{g}(_{0}h_{n})  \; r \; \delta_{x}(_{0}h_{n}) \; d \; \delta_{x}(_{0}h_{n}) \\
_{6}h_{n+1}=\delta_{f} (_{5}h_{n}) \; u \; \overline{ \delta_{m}(_{5}h_{n})} \; r \; _{5}h_{n} \; d \; \overline{ \delta_{y}(_{5}h_{n})} \\
_{7}h_{n+1}=\delta_{f} (_{5}h_{n}) \; u \; \overline{ \delta_{m}(_{5}h_{n})} \; r \; _{5}h_{n} \; d \; \delta_{a} (_{5}h_{n}) \\
_{8}h_{n+1}=\overline{\delta_{g} (_{5}h_{n})} \; u \; \overline{ \delta_{m}(_{5}h_{n})} \; r \; _{5}h_{n} \; d \; \delta_{a}(_{5}h_{n}) \\
_{9}h_{n+1}=\overline{\delta_{o} (_{5}h_{n})} \; u \;  \delta_{g}(_{5}h_{n}) \; r \; \overline{\delta_{a}(_{5}h_{n})} \; d \; \delta_{x} (_{5}h_{n}) \\
_{10}h_{n+1}=\delta_{m} (_{5}h_{n}) \; u \;  \delta_{g}(_{5}h_{n}) \; r \; \overline{\delta_{a}(_{5}h_{n})} \; d \; \overline{_{5}h_{n}} \\
_{11}h_{n+1}=\delta_{m} (_{5}h_{n}) \; u \;  \delta_{g}(_{5}h_{n}) \; r \; \overline{\delta_{a}(_{5}h_{n})} \; d \; \delta_{x}(_{5}h_{n}) 
\end{array}
\end{equation}
 
\section{Locality preservation}\label{sec:locality}

Locality preservation in a HHCK curve, will depend on the curve itself (i.e. the used affine transformations) and our election of the kernel. HHCK curves may show as good locality preservation properties than the HHC ones. To illustrate this, we examine the locality preservation properties of the $HHCK$ built with the mouse and frog kernels.

The dilation factor\cite{bauman06} (also known as square-to-linear ratio) is defined as

\begin{equation}
 \sigma_{f}=\displaystyle \sup_{\begin{array}{c}t,s \in [0,1] \\ t \neq s\end{array}} \frac{\parallel f(t)-f(s)\parallel^2}{|t-s|}
\end{equation}

where $f: \mathbb{R}\rightarrow\mathbb{R}^2$ is the mapping induced by the space filling curve. While $\sigma_{f}$ is strictly defined for the infinite order space filling curve, the same expression can be used to compare $n$th-order approximations of different curves. The dilation factor is closely related to the H\"older continuity condition. It is generally considered that the smaller the dilation factor, the better the local preserving properties of the mapping \cite{shchepin08}.

\begin{figure}
\includegraphics[scale=0.7]{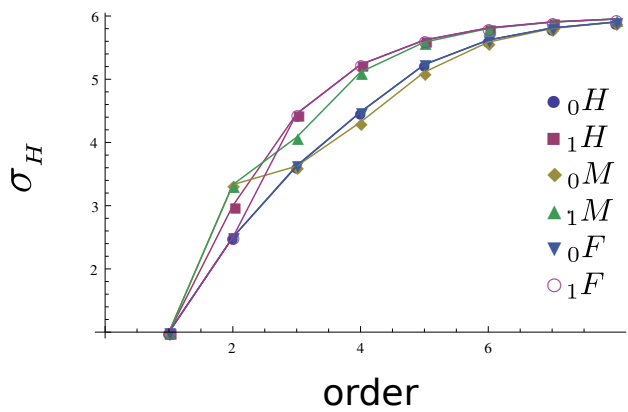}
\caption{Dilation factor with increasing order for the Hilbert curve $_{0}H$, the Moore curve $_{1}H$, the Hilbert curve ($_{0}H_{m}$) and Moore curve ($_{0}H_{m}$) with the mouse kernel, and  the Hilbert curve ($_{0}H_{f}$) and Moore curve ($_{0}H_{f}$) with the frog kernel.
}\label{fig:dilationfactor}
\end{figure}

For the Hilbert curve $_{0}H$, it is known that the dilation factor reaches a value of six in the infinite order limit \cite{bauman06}. The same value is attained for the Moore curve. Figure \ref{fig:dilationfactor} plots for increasing order the dilation factor $\sigma_{H}$ for the Hilbert and Moore curve together with $_{0}M$, $_{1}M$, $_{0}F$ and $_{1}F$. It can be generally seen that the three kernels have similar behavior, with the $_{0}H$, $_{0}M$ and $_{0}F$ curves, growing slightly faster than the $_{1}H$, $_{1}M$ and $_{1}F$ curves, but both sets of curves with very close values for a fixed order.  The introduced two new kernels do not show a degradation of the dilation factor values compared to the HHC.

\begin{figure}
\includegraphics[scale=0.5]{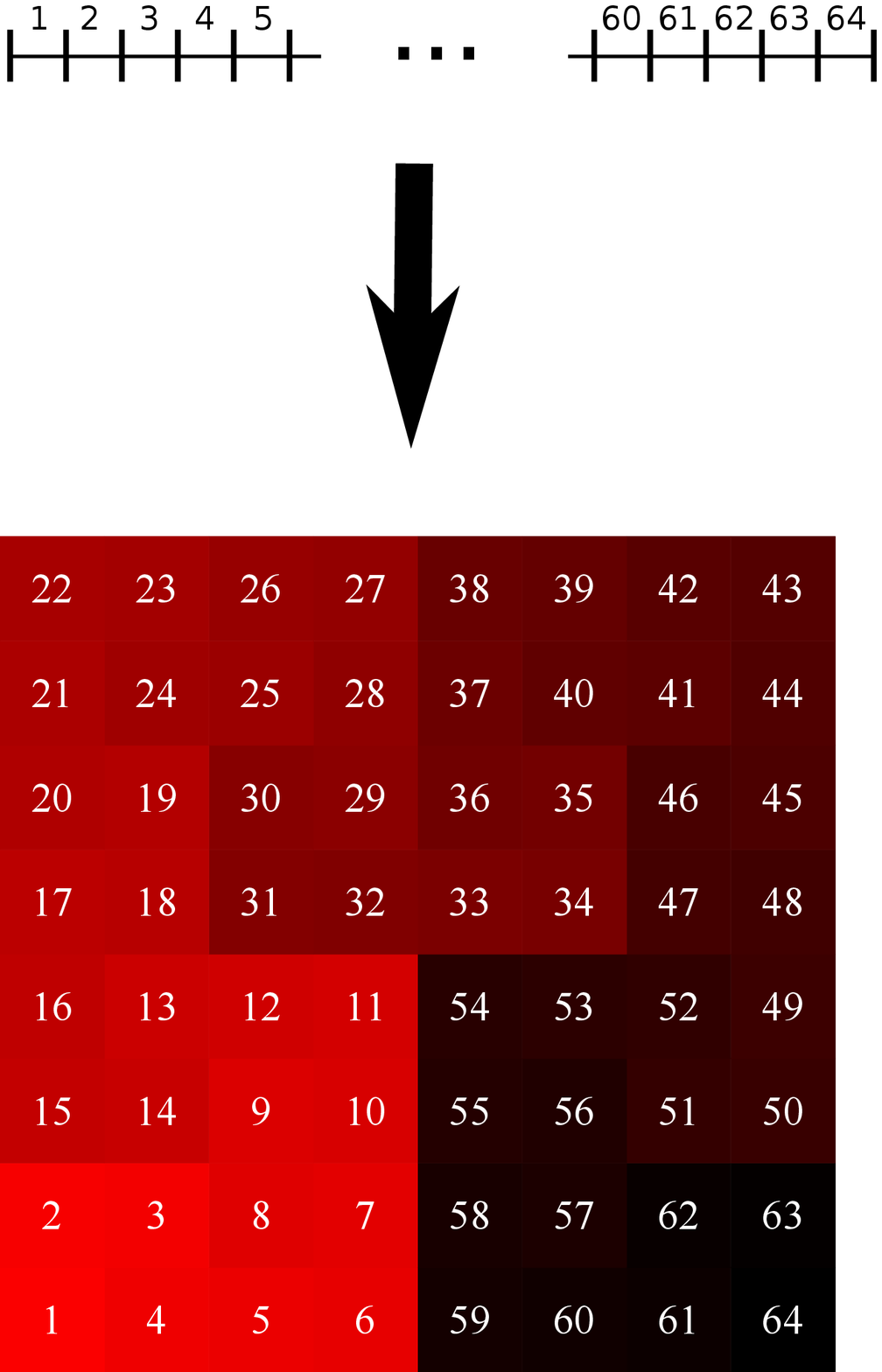}
\caption{The bijective mapping of the $3^{th}$ order iteration of the Hilbert curve $_{0}H_{4}$. The unit interval is partitioned into $4^3=64$ subintervals labeled in increasing order. To each subintervals corresponds a subsquare in the unit square $Q$ which is accordingly labeled.
}\label{fig:hilbertgrid}
\end{figure}

Dilation factor is a global descriptor of locality, it will be important to have a site measure of locality preservation in a space filling curve.  For a curve of order $n$,  we can label each $4^n$ subintervals in increasing integer order, and assign to each of the $4^n$ subsquare under the space filling mapping, the corresponding label (Figure \ref{fig:hilbertgrid}).  Consider at each subsquare the mean value of the difference between its label and those of its eight nearest neighbors. The larger the mean difference value, the less locality is preserved at the site. If we do this calculation for every subsquare, the result is a congruent square to $Q$ where each subsquare, is labeled by the mean difference value, we call such unit square a difference map and denote it by a letter ''d'' in front of the curve symbol. 

The difference map will show, those regions and boundary through which locality is poorly preserved, giving a better picture of the locality preserving properties of the space filling curve.

\begin{figure}
\includegraphics[scale=2]{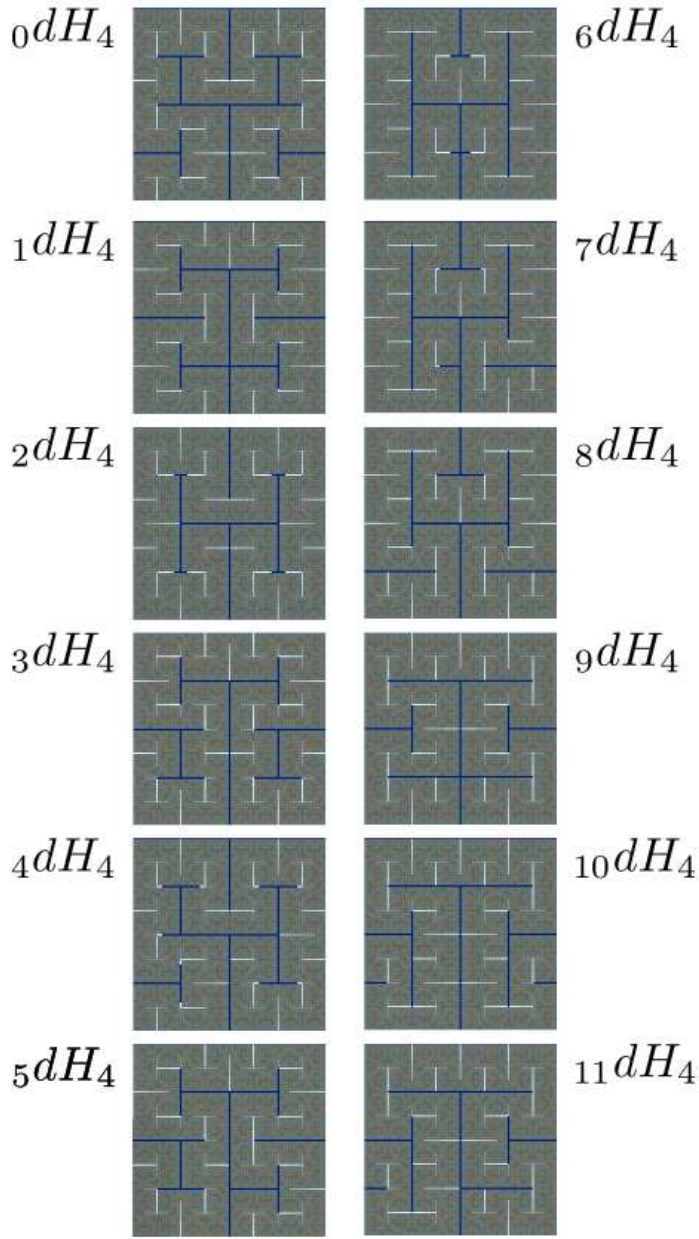}
\caption{The difference map of the 12 HHC. The blue lines corresponds to the largest loss of locality (higher than the mean difference value of the whole curve plus one standard deviation), followed by the white lines. Gray colors are given in a log scale.  
}\label{fig:diffmapHHC}
\end{figure}

Figure \ref{fig:diffmapHHC} shows the difference map for the twelve HHC, similar maps can be drawn for the HHCK  with essentially the same features. We define a strong locality  barrier as the continuous locus of points with values larger than the mean value of the whole map plus the standard deviation (shown in blue in the figure).  The first thing to notice is that difference maps follow the symmetry of the corresponding space filling curve. A strong boundary of locality rupture separating the first and fourth quadrant is observed independent on the type of curve, which should be expected as there is no connectivity between the first and fourth quadrant.  Each difference map gives a good idea on the entry and exit point of each quadrant, for example in the HHC $_{6}H_{4}$ there is a better locality when crossing between the second and third quadrant than in $_{2}H_{4}$. The improper curves $_{6}H_{4}$ to $_{11}H_{4}$, has a better connectivity across the $1\longrightarrow2$ and $3\longrightarrow4$ quadrants, as seen by a shorter (1/2 length for each quadrant boundary) strong locality  barrier, compared to the proper curves, where, except for $_{2}H_{4}$ curve, in all other five cases, the strong locality  barrier is 3/4 of the quadrant boundary length. A similar analysis for the $2\longrightarrow3$ crossover shows that the better connectivity is attained for all improper curves and for the proper curves $_{1}H_{4}$, $_{3}H_{4}$ and $_{5}H_{4}$, where the strong locality  barrier is half the quadrant boundary length, being 3/4 in the the other curves.

\begin{figure}
\includegraphics[scale=0.75]{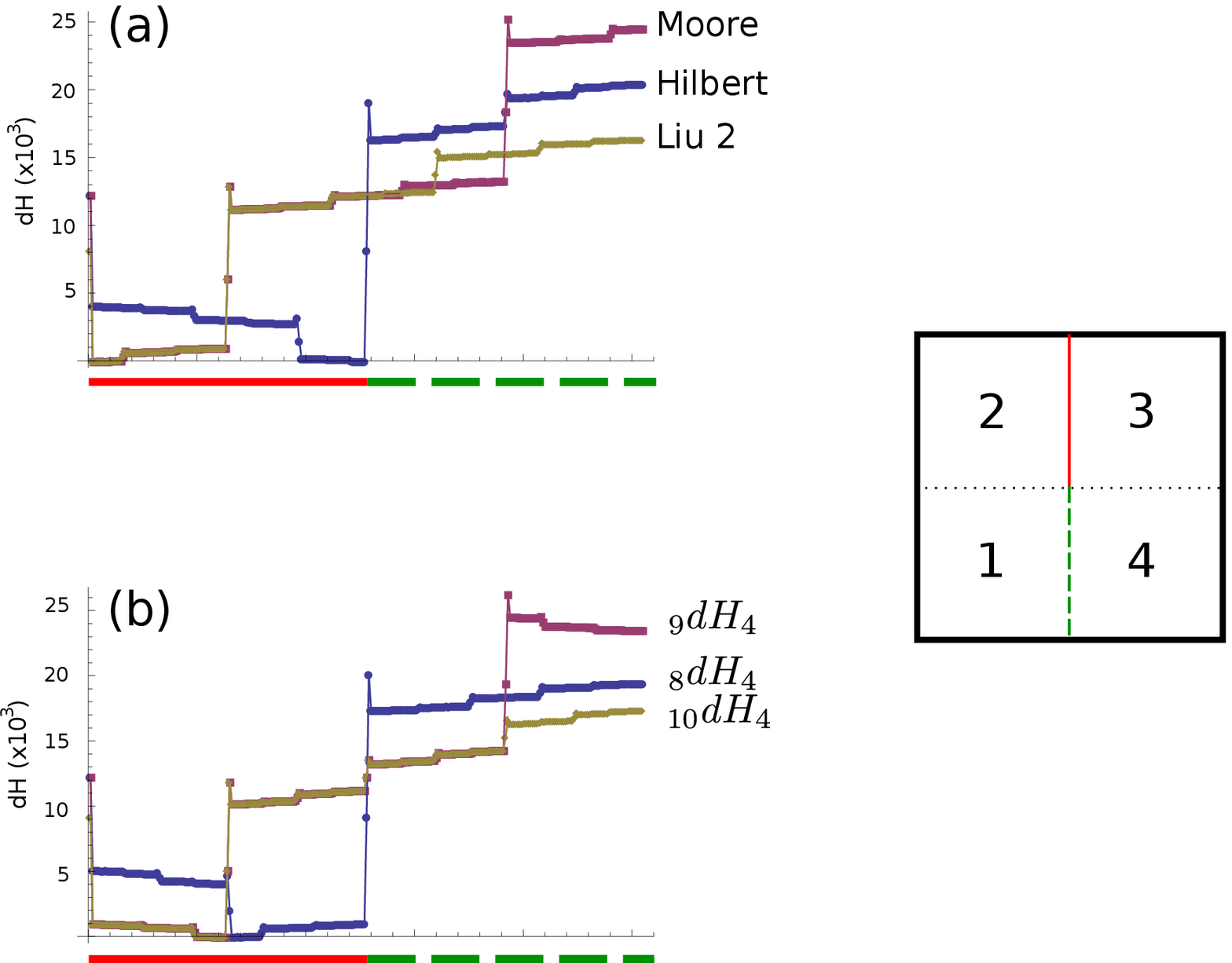}
\caption{The difference map along the boundary between the 1-2 and 3-4 quadrant for (a) the proper HHC $_{0}H_{4}$ (Hilbert), $_{1}H_{4}$ (Moore) and $_{3}H_{4}$ (Liu 2); and the improper HHC (b) $_{8}H_{4}$, $_{9}H_{4}$, $_{10}H_{4}$.
}\label{fig:centralbarrier}
\end{figure}

Comparison of the mean difference values across the $2\longrightarrow3$ and $4\longrightarrow1$ quadrants for several curves is shown in figure \ref{fig:centralbarrier}. In all the analyzed cases, the locality is better preserved along the $2\longrightarrow3$ boundary with mean difference values below one half of those found in the $4\longrightarrow1$ boundary. The better behaving proper curve is $_{3}H_{4}$(Liu 2) with the smaller barrier along both quadrants boundaries, while the worst behaving proper curve is the Moore curve with the highest value locality barrier. Hilbert curve $_{0}H_{4}$ behaves intermediate between the other two proper curves. The locality barrier for all three curves show three abrupt climbs.  While in the Moore and Hilbert curve, the first climb is at half the total boundary length ( where the $4\longrightarrow1$ starts) in the Liu 2 curve, the first jump happens earlier, at 1/4 of the total boundary length. For the Hilbert and Moore curve, the second climb in the locality barrier occurs at half the $4\longrightarrow1$  boundary, the Liu 2 curve jumps earlier at 1/4 of the same boundary. 

For the improper curves, $_{10}H_{4}$ is the better behaved curve and it exhibits four abrupt climbs at half of each quadrants boundaries, $_{8}H_{4}$ and $_{9}H_{4}$, show similar behavior to the Hilbert and Moore curve respectively.

Statistical analysis of the difference map shows that all HHC, and HHCK curves studied in this contribution behaves similarly (Tables \ref{tbl:statHHC}, \ref{tbl:statmouse}, \ref{tbl:statfrog}). For all the $_{\nu}H$, $_{\nu}M$,$_{\nu}F$ curves, the mean value (second column) of the difference map lies between the value of 259 and 283, where the lowest value is attained for the corresponding eighth curve $_{8}H$ ($_{8}M$ and $_{8}F$), and the largest value for the Moore curve ($_{1}H$, $_{1}M$, $_{1}F$). $_{2}H$ ($_{2}M$ and $_{2}F$) shows the worst mean difference value at a site (third column), while $_{3}H$ ($_{3}M$ and $_{3}F$) exhibits the best mean difference value. The lowest value at a site is around 3 for all considered curves (fourth column). For all the curves, the median (fifth column) is around 10, that when compared to the maximum mean difference value at a site, seems to infer that strong locality barriers are scarce in the curves a conclusion that is verified by the fact that 90\% (seventh column) of the difference map values are below the mean value.  It is interesting to note that the median is insensitive to the specific affine transformation defining the curve and changes with the chosen kernel. The mouse kernel curves exhibit a slightly better mean value of 9.87, where in the case of the HHC this values attains it maximum at 10.75. The distribution of mean difference values allows to calculate the Shannon entropy of the difference map. The Shannon entropy will be a measure of homogeneity in the distribution of values, the more uniform the distribution of mean difference values, the larger the Shannon entropy. The mouse kernel curves have the larger entropy values while the HHC shows the lowest entropy values. difference value entropy, similarly to the median, has the same value for all 12 curves of the same kernel.

\section{Conclusions}
Space filling curves based on Hilbert's construction can be further generalized by the use of arbitrary kernels as long as well behaved connectivity is kept. The use of arbitrary kernels can lead to n-order approximations at least as good as the original twelve HHC. The introduction of arbitrary kernels to construct the HHC, further increase the options of scan patterns for different applications. We show such generalization by introducing two new kernels that behave as good as the HHC.  

An important property of space filling curves is its locality properties. Locality measures how well is preserved neighborhood when mapping from the unit interval to the space filling curve. As much as dilation factor has been used as a measure of locality in space filling curves, it is a global descriptor not fitted to  study the site locality properties of a space filling curve. For the latter purpose we introduce the difference map that allows a more detailed study of locality. From the difference map, several global magnitudes can be defined that further enriches the analysis of space filling curves. It would be interesting to compare HHC with other space filling curves, such as Peano construction, with respect to the site locality properties.     

\section{Acknowledgments}

Universidad de la Habana and Universidad de la Ciencias Inform\'aticas are acknowledge for financial support and computational infrastructure. 


\section{References}

\pagebreak

\begin{table}
\caption{The complete set of affine transformations $q=1/2\;[U,\mathbf{t}]$ for the 12 HHC (The overbar means reversion operation, see \cite{estevez_may13} for details).}\centering\small
\begin{tabular}{lllll}
Curve      & $q_{1}$ & $q_{2}$ & $q_{3}$ & $q_{4}$ \\
\hline
\\
Hilbert    & $\left [ U_{R}, t_{0}\right ]$ & $\left [ U_{I}, t_{1}\right ]$ & $\left [ U_{I}, t_{3}\right ]$ & $\left [ -U_{R}, t_{4}\right ]$ \\ 
Moore    & $\left [ U_{V}, t_{2}\right ]$ & $\left [ U_{V}, t_{3}\right ]$ & $\left [ -U_{V}, t_{5}\right ]$ & $\left [ -U_{V}, t_{3}\right ]$ \\
Liu 1    & $\left [ -U_{I}, t_{3}\right ]$ & $\left [ U_{I}, t_{1}\right ]$ & $\left [ U_{I}, t_{3}\right ]$ & $\left [ -U_{I}, t_{4}\right ]$ \\
Liu 2    & $\left [ U_{H}, t_{1}\right ]$ & $\left [ U_{V}, t_{3}\right ]$ & $\left [ -U_{V}, t_{5}\right ]$ & $\left [U_{H}, t_{3}\right ]$ \\
Liu 3    & $\left [ U_{R}, t_{0}\right ]$ & $\left [ U_{I}, t_{1}\right ]$ & $\left [ U_{I}, t_{3}\right ]$ & $\left [-U_{I}, t_{4}\right ]$ \\
Liu 4    & $\left [ U_{H}, t_{1}\right ]$ & $\left [ U_{V}, t_{3}\right ]$ & $\left [ -U_{V}, t_{5}\right ]$ & $\left [-U_{V}, t_{3}\right ]$ \\
Improper 1    & $\left [ -U_{I}, t_{3}\right ]$ & $\overline{\left [ -U_{H}, t_{3}\right ]}$ & $\left [ U_{I}, t_{3}\right ]$ & $\overline{\left [U_{H}, t_{3}\right ]}$ \\
Improper 2    & $\left [ -U_{I}, t_{3}\right ]$ & $\overline{\left [ -U_{H}, t_{3}\right ]}$ & $\left [ U_{I}, t_{3}\right ]$ & $\left [-U_{R}, t_{4}\right ]$ \\
Improper 3    & $\overline{\left [ -U_{V}, t_{1}\right ]}$ & $\overline{\left [ -U_{H}, t_{3}\right ]}$ & $\left [ U_{I}, t_{3}\right ]$ & $\left [-U_{R}, t_{4}\right ]$ \\
Improper 4    & $\overline{\left [ -U_{R}, t_{3}\right ]}$ & $\left [ U_{V}, t_{3}\right ]$ & $\overline{\left [ U_{R}, t_{3}\right ]}$ & $\left [-U_{V}, t_{3}\right ]$ \\
Improper 5    & $\left [ U_{H}, t_{1}\right ]$ & $\left [ U_{V}, t_{3}\right ]$ & $\overline{\left [ U_{R}, t_{3}\right ]}$ & $\overline{\left [-U_{I}, t_{4}\right ]}$ \\
Improper 6    & $\left [ U_{H}, t_{1}\right ]$ & $\left [ U_{V}, t_{3}\right ]$ & $\overline{\left [ U_{R}, t_{3}\right ]}$ & $\left [-U_{V}, t_{3}\right ]$ \\
\\
\hline
\end{tabular}\label{tbl:affine}

\begin{tabular}{l}
\\
\end{tabular}

\begin{tabular}{ll}
$U_{I}=\left ( \begin{array}{cc} 1 & 0  \\ 0 & 1 \end{array} \right )$ &
$U_{R}=\left ( \begin{array}{cc} 0 & 1  \\ 1 & 0 \end{array} \right ) $ \\
$U_{V}=\left ( \begin{array}{cc} 0 & -1  \\ 1 & 0 \end{array} \right )$ &
$U_{H}=\left ( \begin{array}{cc} 1 & 0  \\ 0 & -1 \end{array} \right )$
\end{tabular}

\begin{tabular}{l}
\\
\end{tabular}

\begin{tabular}{lll}
$t_{0}=\left ( \begin{array}{c} 0 \\ 0 \end{array} \right ) $ &
$t_{1}=\left ( \begin{array}{c} 0 \\ 1 \end{array} \right ) $ &
$t_{2}=\left ( \begin{array}{c} 1 \\ 0 \end{array} \right ) $ \\
$t_{3}=\left ( \begin{array}{c} 1 \\ 1 \end{array} \right ) $ &
$t_{4}=\left ( \begin{array}{c} 2 \\ 1 \end{array} \right ) $ &
$t_{5}=\left ( \begin{array}{c} 1 \\ 2 \end{array} \right ) $ 
\end{tabular}
\end{table}

\begin{table}
\caption{Statistical properties of the difference map for the HHC (See text for details).}
\begin{tabular}{llllllll}
Curve & Mean & Max & Min & Median & Entr. & \% $<$ mean \\
\hline \\
$_{0}dH_{4}$  & 262 & 20480 & 3 & 10.75 & 5.71  & 90.0   \\
$_{1}dH_{4}$  & 269 & 25259 & 3 & 10.75 & 5.72  & 90.2   \\
$_{2}dH_{4}$  & 283 & 28672 & 3 & 10.75 & 5.71  & 90.8   \\
$_{3}dH_{4}$  & 262 & 16384 & 3 & 10.75 & 5.71  &  89.9  \\
$_{4}dH_{4}$  & 272 & 25941 & 3 & 10.75 & 5.72  &  90.8  \\
$_{5}dH_{4}$  & 265 & 22870 & 3 & 10.75 & 5.73  &  90.2  \\
$_{6}dH_{4}$  & 280 & 27307 & 3 & 10.75 & 5.73  &  90.9  \\
$_{7}dH_{4}$  & 269 & 25259 & 3 & 10.75 & 5.74  &  90.3  \\
$_{8}dH_{4}$  & 259 & 20139 & 3 & 10.75 & 5.73  &  89.9  \\
$_{9}dH_{4}$  & 267 & 26282 & 3 & 10.75 & 5.73  &  90.3  \\
$_{10}dH_{4}$ & 259 & 17408 & 3 & 10.75 & 5.73  &  89.9  \\
$_{11}dH_{4}$ & 263 & 23552 & 3 & 10.75 & 5.73  &  90.1  \\
\end{tabular}\label{tbl:statHHC}
\end{table}

\begin{table}
\caption{Statistical properties of the difference map for the mouse kernel HHCK (See text for details).}
\begin{tabular}{llllllll}
Curve & Mean & Max & Min & Median & Entr. & \% $<$ mean \\
\hline \\
$_{0}dM_{4}$  & 262  & 20480  & 3.25 & 9.87  & 6.48  & 90.0  \\
$_{1}dM_{4}$  & 269  & 25258  & 3.25 & 9.87  & 6.48  & 90.2  \\
$_{2}dM_{4}$  & 283  & 28672  & 3.25 & 9.87  & 6.48  & 90.8  \\
$_{3}dM_{4}$  & 262  & 16384  & 3.25 & 9.87  & 6.48  & 89.9  \\
$_{4}dM_{4}$  & 272  & 25941  & 3.25 & 9.87  & 6.48  & 90.8  \\
$_{5}dM_{4}$  & 265  & 22869  & 3.25 & 9.87  & 6.49  & 90.1  \\
$_{6}dM_{4}$  & 280  & 27306  & 3.00 & 9.87  & 6.50  & 90.9  \\
$_{7}dM_{4}$  & 269  & 25258  & 3.00 & 9.87  & 6.50  & 90.3  \\
$_{8}dM_{4}$  & 259  & 20139  & 3.25 & 9.87  & 6.48  & 89.9  \\
$_{9}dM_{4}$  & 267  & 26283  & 3.00 & 9.87  & 6.49  & 90.3  \\
$_{10}dM_{4}$ & 259  & 17407  & 3.00 & 9.87  & 6.48  & 89.9  \\
$_{11}dM_{4}$ & 263  & 23552  & 3.00 & 9.87  & 6.50  & 90.1  \\
\end{tabular}\label{tbl:statmouse}
\end{table}

\begin{table}
\caption{Statistical properties of the difference map for the frog kernel HHCK (See text for details).}
\begin{tabular}{llllllll}
Curve & Mean & Max & Min & Median & Entr. & \% $<$ mean \\
\hline \\
$_{0}dF_{4}$  & 262  & 20480  & 2.75 & 10.5 & 6.18  & 90.0  \\
$_{1}dF_{4}$  & 269  & 25259  & 2.75 & 10.5 & 6.18  & 90.2  \\
$_{2}dF_{4}$  & 283  & 28672  & 2.75 & 10.5 & 6.18  & 90.8  \\
$_{3}dF_{4}$  & 262  & 16384  & 2.75 & 10.5 & 6.17  & 89.9  \\
$_{4}dF_{4}$  & 272  & 25942  & 2.75 & 10.5 & 6.18  & 90.8  \\
$_{5}dF_{4}$  & 265  & 22870  & 2.75 & 10.5 & 6.19  & 90.2  \\
$_{6}dF_{4}$  & 279  & 27307  & 2.75 & 10.5 & 6.18  & 90.9  \\
$_{7}dF_{4}$  & 269  & 25259  & 2.75 & 10.5 & 6.19  & 90.3  \\
$_{8}dF_{4}$  & 259  & 20139  & 2.75 & 10.5 & 6.18  & 89.9  \\
$_{9}dF_{4}$  & 267  & 26283  & 2.75 & 10.5 & 6.18  & 90.3  \\
$_{10}dF_{4}$ & 259  & 17408  & 2.75 & 10.5 & 6.18  & 89.9  \\
$_{11}dF_{4}$ & 263  & 23553  & 2.75 & 10.5 & 6.19  & 90.1  \\
\end{tabular}\label{tbl:statfrog}
\end{table}


\end{document}